\magnification=1200

\input amstex

\documentstyle{amsppt}

%%%%%%%%%%%%%%%%%%%%%%%%%%%%%%%%%%%%%%%%%%%%%%%%%%%%%%%%%%%%%%
\pagewidth{165truemm}
\pageheight{227truemm}

%%%%%%%%%%%%%%%%%%%%%%%%%%%%%%%%%%%%%%%%%%%%%%%%%%%%%%%%%%%%%%

\def\C{{{\Bbb C}}}

\def\p#1{{{\Bbb P}^{#1}_{k}}}

\def\dualp#1{{{\check{\Bbb P}}^{#1}_{k}}}

\def\Pic{\operatorname{Pic}}

\def\Ext{\operatorname{Ext}}

\def\Sing{\operatorname{Sing}}

\def\GL{\operatorname{GL}}

\def\im{\operatorname{im}}

\def\char{\operatorname{char}}

\def\ap{\operatorname{ap}}

\def\Ofa#1{{{\Cal O}_{#1}}}

\def\Hom#1{{{\Cal H}\kern -0.25ex{\italic om\/}_{\Ofa#1}}}

\def\mapright#1{\mathbin{\smash{\mathop{\longrightarrow}\limits^{#1}}}}

\def\ga#1{{\accent"12 #1}}

\def\Hilb{{{\Cal H}\kern -0.25ex{\italic ilb\/}}}

\topmatter

\title
Canonical curves with low apolarity
\endtitle

\rightheadtext{}

\author 
E. Ballico, G. Casnati, R. Notari
\endauthor

\address
Edoardo Ballico, Dipartimento di Matematica, Universit\`a degli Studi di Trento,
Via Sommarive 14, 38123 Povo, Italy
\endaddress

\email
ballico\@science.unitn.it
\endemail

\address Gianfranco Casnati, Dipartimento di Matematica, Politecnico di Torino,
c.so Duca degli Abruzzi 24, 10129 Torino, Italy
\endaddress

\email
casnati\@calvino.polito.it
\endemail

\address
Roberto Notari, Dipartimento di Matematica \lq\lq Francesco Brioschi\rq\rq, Politecnico di Milano,
via Bonardi 9, 20133 Milano, Italy
\endaddress

\email
roberto.notari\@polimi.it
\endemail



\abstract 
Let $k$ be an algebraically closed field and let $C$ be a non--hyperelliptic smooth projective curve of genus $g$ defined over $k$. Since the canonical model of $C$ is arithmetically Gorenstein, Macaulay's theory of inverse systems allows to associate
to $C$ a cubic form $f$ in the divided power $k$--algebra $R^{g-3}$ in $g-2$ variables. The apolarity $\ap(C)$ of $C$ is
the minimal number $t$ of linear form $\ell_1,\dots,\ell_t\in R^{g-3}$ needed to write $f$ as sum of their divided power cubes. 

It is easy to see that  $\ap(C)\ge g-2$ and P\. De Poi and F\. Zucconi classified curves with $\ap(C)=g-2$ when $k\cong\C$. In this paper, we give a complete, characteristic free, classification of curves $C$  with apolarity $g-1$ (and $g-2$). 
\endabstract

\endtopmatter
\document

\head
0. Introduction and notation
\endhead
Throughout the paper $k$ will denote an algebraically closed field of arbitrary characteristic. Let $C$ be a non--hyperelliptic smooth projective curve of genus $g\ge3$ over $k$. It is well--known that the canonical map $\phi\colon C\to\vert\omega_C\vert\check{\ }=\dualp{g-1}$ is an embedding. 

The classical Babbage--Enriques--Noether--Petri Theorem (see [SD]) gives a complete description of the generators of its homogeneous ideal $I_{\phi(C)}\subseteq k[x_0,\dots,x_{g-1}]$ and some results are known about the higher syzygies of its homogeneous coordinate ring $T(\phi(C)):=k[x_0,\dots,x_{g-1}]/I_{\phi(C)}$, expecially in the set up of the well--known Green's conjecture.

Another approach for studying the canonical embedding of $C$ is suggested in [I--R] (see also [R--S]) when $k\cong\C$ is the complex field. The basic idea is that the ring $T(\phi(C))$ is Gorenstein, i.e. it is self--injective. Since $T(\phi(C))$ has Krull--dimension $2$, it follows that for any choice of general linear forms $h_1,h_2\in k[x_0,\dots,x_{g-1}]$, the quotient ring $T(\phi(C))/(h_1,h_2)$ is an Artinian Gorenstein graded ring and it is not difficult to check that its Hilbert function is $(1,g-2,g-2,1)$.  

Artinian Gorenstein graded rings with Hilbert function $(1,n,n,1)$ have been object of a deep study (see e.g. [Ia1] or [C--R--V] for general results about them, [C--N1] for their classification when $n\le3$ and [C--N2] for some results in the case $n=4$). In particular the well--known Macaulay's inverse system method asserts that there is a bijective correspondence between Artinian Gorenstein graded quotient rings of $k[t_1,\dots,t_n]$ with Hilbert function $(1,n,n,1)$ and forms of degree $3$ up to scalars in the divided power $k$--algebra in the variables $y_1,\dots,y_n$ (see Appendix A of [I--K] for the definition and main properties of the divided power algebra: for the sake of simplicity we only recall that such an algebra is isomorphic to $k[y_1,\dots,y_n]$ when $k\cong\C$).

It is thus natural to relate geometric properties of $C$, hence projective properties of its canonical model $\phi(C)$, with algebraic properties of the corresponding quotient ring $T(\phi(C))/(h_1,h_2)$. When $k\cong \C$, such an idea has been smartly used in [DP--Z1] in order to prove the following result (see Theorem 4: see also [DP--Z2] for similar results for subcanonical curves).

\proclaim{Theorem}
Let $k\cong \C$. A non--hyperelliptic  smooth projective curve $C$ of genus $g\ge3$ over $k$ is either trigonal or isomorphic to a plane quintic if and only if the corresponding polynomial is a Fermat cubic in $k[y_0,\dots,y_{g-3}]$ up to the natural action of the general linear group $\GL_{g-2}$ (i.e. the corresponding polynomial is in the $\GL_{g-2}$--orbit of $y_0^3+\dots+y_{g-3}^3$).
\qed
\endproclaim

The proof of the above Theorem rests on the above mentioned description of the generators of the ideal $I_{\phi(C)}$. Since such a description is characteristic free our first result is its generalization in any characteristic. In particular in Section 1 we prove that 

\proclaim{Theorem A}
A non--hyperelliptic  smooth projective curve $C$ of genus $g\ge3$ over $k$ is either trigonal or isomorphic to a plane quintic if and only if the corresponding polynomial is a Fermat cubic in the divided power algebra up to the natural action of the general linear group $\GL_{g-2}$.
\qed
\endproclaim

It is interesting to recall that the above mentioned Babbage--Enriques--Noether--Petri Theorem relates the degrees of a minimal set of generators of $I_{\phi(C)}$ with the existence of surfaces $S\subseteq\dualp{g-1}$ of minimal degree $g-2$ containing $\phi(C)$. 

Thus it is natural to inspect the case when the canonical model of $C$ lies on a surface $S\subseteq\dualp{g-1}$ of almost minimal degree $g-1$. Such an analysis yields to the main result of this paper proved in Section 4.

\proclaim{Theorem B}
A non--hyperelliptic  smooth projective curve $C$ of genus $g\ge5$ over $k$ is either bielliptic (i.e. there is a map of degree $2$ from $C$ to a smooth elliptic curve) or it is birationally isomorphic to a plane sextic with at most $10-g$ double points as singularities (whose position satisfies an extra technical condition)  or it is isomorphic to a smooth complete intersection inside $\p3$ of an integral quadric $Q$ with an integral quartic surface if and only if the corresponding polynomial can be written as the sum of $g-1$ cubes in the divided power algebra.
\qed
\endproclaim

The structure of this paper is as follows. In Section 1 we recall some facts about apolarity and we prove Theorem A above. In Section 2 we summarize some results about surfaces of almost minimal degree: in particular we give their classification in the normal case. In Section 3 we prove a characterization of canonical curves lying on a surface of almost minimal degree. Finally, in Section 4, we prove the main result of this paper, namely Theorem B above.

The authors would like to thank E\. Carlini, F.O\. Schreyer and F\. Vaccarino for their helpful comments and fruitful suggestions about the contents and the exposition of this paper.

\subhead 
Notation
\endsubhead
We work over an algebraically closed field $k$ and we will denote by $\char(k)$ its characteristic. The symbol $\GL_n$  denotes the general linear group of invertible $n\times n$ matrices with coefficients in $k$.

Let $k[t_0,\dots,t_n]$ be the ring of polynomials in the variables $t_0,\dots,t_n$ with coefficients in $k$ an let $k[t_0,\dots,t_n]_d$ be the
vector space of degree $d$ forms. If $q_1,\dots,q_t\in k[t_0,\dots,t_n]_d$, we denote by $\langle q_1,\dots,q_t\rangle$ the subspace of $k[t_0,\dots,t_n]_d$ generated by $q_1,\dots,q_t$. 

If $V$ is a vector space, then we denote by ${\Bbb P}(V)$ the corresponding projective space. In particular we set $\p n:={\Bbb P}(k^{n+1})$. A curve is a projective connected smooth scheme of dimension $1$. 

Let $\p n$ be the projective space with coordinates $t_0,\dots,t_n$. If $q_1,\dots,q_t\in k[t_0,\dots,t_n]$ are forms we will denote by $D_0(q_1,\dots,q_t)$ the corresponding zero scheme inside $\p n$. If $X\subseteq\p n$ is a closed subscheme, then we denote by $I_X$ its saturated homogeneous ideal in $k[t_0,\dots,t_n]$ and we define its homogeneous coordinate ring as 
$$
T(X):=k[t_0,\dots,t_n]/I_X.
$$
We will say that $X$ is arithmetically Cohen--Macaulay (resp. arithmetically Gorenstein), aCM for short (resp. aG for short), in $\p n$ if the ring $T(X)$ is Cohen--Macaulay (resp. Gorenstein).

If $\gamma:=(\gamma_0,\dots,\gamma_n)\in{\Bbb N}^{n+1}$ is a multi--index, then we set
$\vert \gamma\vert:=\sum_{i=0}^n\gamma_i$, $\gamma!:=\prod_{i=0}^n{\gamma_i!}$, $t^\gamma:=t_0^{\gamma_0}\dots t_n^{\gamma_n}\in k[t_0,\dots,t_n]$ and we say that $\gamma\ge0$ if and only if $\gamma_i\ge0$ for each $i=0,\dots,n$: if $\delta:=(\delta_0,\dots,\delta_n)\in{\Bbb N}^{n+1}$ is another multi--index then we write $\gamma\ge\delta$ if and only if $\gamma-\delta\ge0$.

For other definitions, results and notation we always refer to [Ha].

\head
1. Apolarity and first results
\endhead
In [I--K] (see also [R--S] and [I--R] for the characteristic $0$ case) some facts about the classical Macaulay correspondence are summarized. 

We recall that we can consider two graded rings, namely the polynomial ring  $T^n:=k[x_0,\dots,x_n]$ and the divided power $k$--algebra $R^n$ in the $n+1$ variables $y_0,\dots,y_n$. 

As explained in [I--K], Appendix A, the $k$--vector space $R^n$ coincides with $k[y_0,\dots,y_n]$, hence there exists a natural action of $\GL_{n+1}$ on $R^n$. Let $\gamma:=(\gamma_0,\dots,\gamma_n)\in{\Bbb N}^{n+1}$: the divided power monomial $y_0^{\gamma_0}\dots y_n^{\gamma_n}\in R^n$ will be usually denoted as $y^{[\gamma]}$ (instead of $y^{\gamma}$). The $k$--algebra structure on $R^n$ is obtained by extending by linearity the monomial multiplication 
$$
y^{[\gamma]}y^{[\delta]}:={{(\gamma+\delta)!}\over{\gamma!\delta!}}y^{[\gamma+\delta]}.
$$

The algebra $T^n$ acts on $R^n$ by differentiation. More precisely, if $\gamma:=(\gamma_0,\dots,\gamma_n),\delta:=(\delta_0,\dots,\delta_n)\in{\Bbb N}^{n+1}$ we have the natural pairing $T^n_g\times R^n_d\to R^n_{d-g}$ (we denote by $*^n_d$ the summand of degree $d$ elements in $*^n$) given on monomials by the rule
$$
x^{\gamma}(y^{[\delta]}):=\cases y^{[\delta-\gamma]}\qquad&\text{if $\delta\ge\gamma$,}\\
0\qquad&\text{if $\delta\not\ge\gamma$.}
\endcases
$$

For any linear form $\ell:=\sum_{i=0}^na_iy_i\in R^n_1$ the divided power $\ell^{[d]}$ is defined as 
$$
\ell^{[d]}:=\sum_{\delta\in{\Bbb N}^{n+1},\ \vert\delta\vert=d}a_0^{\delta_0}\dots a_n^{\delta_n}y_0^{[\delta_0]}\dots y_n^{[\delta_n]}.
$$
Let $\ell\in R^n_1$ be as above and let us denote by $a:=(a_0,\dots,a_n)$. If $q\in T^n_e$ then the above formula yields $q(\ell^{[d]})=q(a)\ell^{[d-e]}$ if $e\le d$, thus we deduce 
$$
q(\ell^{[d]})=0\iff q(a)=0.\tag 1.1
$$

Such an action defines a perfect pairing $T^n_d\times R^n_d\to k$ between forms of degree $d$ in $R^n$ and in $T^n$. In particular $R^n_1$ and $T^n_1$  are natural dual vector spaces. Therefore the projective spaces with coordinates $y_0,\dots,y_n$ and $x_0,\dots,x_n$ are naturally dual each other and we denote them by $\p n$ and $\dualp n$ respectively.

We will say that {two \sl homogeneous forms $f\in R^n$ and $q\in T^n$ are apolar if $q(f)=0$}\/. As explained in [I--K], apolarity allows us to associate an Artinian Gorenstein graded quotient of $T^n$ to a form in $R^n$ as follows. Let $f\in R^n_d$: then we set
$$
f^\perp:=\{\ q\in T^n\ \vert\ q(f) =0\ \}
$$
and it is easy to prove that both $f^\perp$ is a homogeneous ideal in $T^n$ and $A^f:=T^n/f^\perp$ is an Artinian Gorenstein graded quotient of $T^n$ with socle $0:_{T^n}T^n_1$ in degree $d$. Also the converse is true i.e. if $A$ is an Artinian Gorenstein graded quotient of $T^n$, say $A:=T^n/I$, with socle in degree $d$ then there exists $f\in R^n_d$ such that $I=f^\perp$. The main result about apolarity due to Macaulay (see [I--K], Lemma 2.12 and the references cited there) is the following

\proclaim{Theorem 1.2}
The map $f\mapsto A^f$ induces a bijection between ${\Bbb P}(R^n_d)$ and the set of graded Artinian Gorenstein quotient rings of $T^n$ with socle in degree $d$.
\qed
\endproclaim

Now we restrict our attention to non--hyperelliptic curves $C$ (recall that a curve is smooth by definition in this paper) of genus $g\ge3$. Then the homogeneous coordinate ring of the canonical model of $C$ inside $\dualp{g-1}$ satisfies
$$
T(\phi(C))\cong T_C:=\bigoplus_{h=0}^{+\infty} H^0\big(C,\omega_C^h\big).
$$
If we take a general linear subspace $ H = D_0(h_1, h_2) \subseteq\dualp{g-1}$ of codimension $ 2 $, then the algebra $ T_C / (h_1, h_2) $ turns out to be naturally an Artinian Gorenstein graded quotient of $ T^{g-1} / (h_1, h_2) $ with socle in degree $ 3 $ and Hilbert function $ (1, g-2, g-2, 1)$. With a proper choice of coordinates $ y_0, \dots, y_{g-1} $ in $\p{g-1}$ and of the corresponding dual coordinates $ x_0, \dots, x_{g-1} $ in $\dualp{g-1}$, we can assume that  $h_i=x _{g-3+i}$, $i=1,2$, i.e.
 $ H = V(x_{g-2}, x_{g-1}) $ so that the morphism $T^{g-1}\to T^{g-3}$ sending to zero $x_{g-2}$ and $x_{g-1}$ induces a natural identification $ T^{g-1} / (h_1, h_2)= T^{g-3}$.

Thus we can associate a polynomial $f\in R^{g-3}_3$ to $C$ and $H$ in such a way that, via the aformentioned identification,
$$
f^\perp={{I_{\phi(C)}+(h_1, h_2)}\over{(h_1, h_2)}}
$$
or, equivalently,  
$$
{{T_C}\over{(h_1, h_2)}}={{T^{g-3}}\over{f^\perp}}.
$$
\remark{Remark 1.3}
Since $x_{g-2},x_{g-1}$ is a regular sequence in $T^{g-1}$, it follows that the Betti numbers of $T_C$ as module over $T^{g-1}$
coincide with the Betti numbers of $T_C/(x_{g-2},x_{g-1})\cong T^{g-3}/f^\perp$ as module over $T^{g-3}$ (see e.g. Theorem 1.3.6 of [Mi]).
\endremark
\medbreak

\definition{Definition 1.4}
With the notation introduced above we say that the polynomial $f\in R^{g-3}$ and the curve $C$ are apolar each other.
\enddefinition

Following the definition given in [I--R] when $\char(k)=0$, apolarity allows us to define a rational map 
$$
\psi_C \colon G(g-2,g)\dashrightarrow H_{g-3,3}
$$
where $ G(g-2,g)$ denotes the Grassmannian of subspaces of codimension $2$ of $\dualp{g-1}$ and $H_{g-3,3}:=R^{g-3}_3/\GL_{g-2}$ without restrictions on $\char(k)$. On the other hand in $H_{g-3,3}$ it is also defined the locus $H_{g-3,3}(h)$ of $\GL_{g-2}$--orbits of cubics (in the sense of divided power algebra when $\char(k)>0$) which can be written as sum of $h$ (divided power) cubes of linear forms. Obviously we must have $h\ge1$ and $H_{g-3,3}(h)\subseteq H_{g-3,3}(h+1)$. 
For simplicity we denote by $\ap(C)$ the {\sl apolarity of $C$}\/, i.e. the smallest integer $h$ such that the general point of $\im(\psi_C)$ is in $H_{g-3,3}(h)$. 

When $\char(k)\ne2,3$, Theorem 2 of [Ia2] (which is based on the results proved in [A--H]) asserts that the general point of $H_{g-3,3}$ lies in $H_{g-3,3}(h_{gen})$ where
$$
h_{gen}:=\left\lceil{{g(g-1)}\over6}\right\rceil
$$
except for $g=7$, when $h_{gen}=8$ (instead of $7$). Thus it is reasonable to hope that $\ap(C)\le h_{gen}$. Actually a better upper bound holds: indeed in Remark 3.3 of [I--R] it is proved that $\ap(C)\le2g-4$ when $C$ is general and $\char(k)=0$. 

The following lemma provides a lower bound for the apolarity of a non--hyperelliptic curve $C$.

\proclaim{Lemma 1.5}
Let $C$ be a non--hyperelliptic curve of genus $g\ge3$ apolar to the polynomial $f\in R^{g-3}_3$. If there exist  $\ell_0,\dots,\ell_s\in R^{g-3}_1$ such that $f=\sum_{i=0}^s\ell_i^{[3]}$, then $\dim\langle \ell_0,\dots,\ell_s\rangle= g-2$.
\endproclaim
\demo{Proof}
We have  $\dim\langle \ell_0,\dots,\ell_s\rangle\le g-2$. Assume that $\dim\langle \ell_0,\dots,\ell_s\rangle=t\le g-3$. Up to $\GL_{g-2}$ we can assume that $f\in R^t\subseteq R^{g-3}$, so that $y_{g-3}$ does not appear in $f$. It follows that $ x_{g-3}(f)=0$, thus $x_{g-3}\in f^\perp$. Due to Remark 1.3  there should exist a linear form in $I_{\phi(C)}\subseteq T^{g-1}$, a contradiction. 
\qed
\enddemo

A natural problem is to  inspect curves $C$ with low apolarity, i.e. such that the general point of $\im(\psi_C)$ is in $H_{g-3,3}(k+1)\setminus H_{g-3,3}(k)$. To this purpose a  fundamental tool is the following fundamental result known as {\sl Apolarity Lemma}\/.

\proclaim{Lemma 1.6}
Let $\ell_i\in R^n_1$ and let $L_i\in\dualp n$ be the corresponding point, $i=1,\dots,s$. Then $f=\sum_{i=0}^s\lambda_i\ell_i^{[d]}$ for some $\lambda_i\in k^*$, $i=0,\dots,s$, if and only if the homogeneous ideal $I_\Gamma\subseteq R^n$ of $\Gamma:=\{\ L_0,\dots,L_s\ \}\subset\dualp{n}$ satisfies $I_\Gamma\subseteq f^\perp$.
\endproclaim
\demo{Proof}
In [I--K], Lemma 1.15 it is proved that $\cap_{i=0}^s(\ell_i^{[d]})^\perp=I_\Gamma\cap T^n_d\subseteq T^n_d$. If $I_\Gamma\subseteq f^\perp$, then $I_\Gamma\cap T^n_d\subseteq f^\perp\cap T^n_d$: since the map $T^n_d\times R^n_d\to k$ is a perfect pairing then $f\in \langle \ell_0^{[d]},\dots,\ell_s^{[d]}\rangle$, thus $f=\sum_{i=0}^s\lambda_i\ell_i^{[d]}$ for some $\lambda_i\in k^*$, $i=0,\dots,s$.

Conversely let $f=\sum_{i=0}^s\lambda_i\ell_i^{[d]}$, for some $\lambda_i\in k^*$, $i=0,\dots,s$. If $q\in I_\Gamma$, then $q(\ell_i^{[d]})=0$ due to Formula (1.1), thus $q\in f^{\perp}$.
\qed
\enddemo

The condition $I_\Gamma\subseteq f^\perp$ is often summarized in the standard literature by saying that {\sl $\Gamma$ is apolar to $f$}\/.

The above result and Babbage--Enriques--Noether--Petri Theorem describing the ideal of a canonical curve have been used in [DP--Z1] (see Theorem 4) to prove an interesting characterization of curves which are either trigonal or isomorphic to plane quintic. Here we rephrase it in terms of apolarity of $C$ and we prove it without any restriction on the characteristic of the base field $k$ (see Theorem A stated in the introduction).

\proclaim{Theorem 1.7}
A non--hyperelliptic curve $C$ of genus $g\ge3$ is either trigonal or it is isomorphic to a plane quintic if and only if $\ap(C)=g-2$.
\endproclaim
\demo{Proof}
If $g=3$ the statement is obvious. Indeed each non--hyperelliptic curve $C$ of genus $g$ is trigonal. On the other hand the ring $T_C/(x_{1},x_{2})$ has Hilbert function $(1,1,1,1)$, thus
$$
T_C/(x_1,x_2)\cong T^{0}/(x_0^4)\cong T^{0}/f^\perp
$$
where $f:=y_0^{[3]}$.

Thus, from now on, we will assume $g\ge4$. It is easy to check that $f$ is the, possibly divided power, Fermat cubic $y_0^{[3]}+\dots+y_{g-3}^{[3]}\in R^{g-3}_3$ if and only if 
$$
f^{\perp} =(x_i x_j , x_i^3-x_{g-3}^3)_{0 \le i < j \le g-3}.
$$
Due to Remark 1.3 this happens if and only if a minimal set of generators of the ideal $I_{\phi(C)}$ contains at least one cubic form. The well--known Babbage--Enriques--Noether--Petri Theorem on the ideal of the canonical model of $C$, which holds true in any characteristic for curves (see [SD], Theorem 3.1), says that $I_{\phi(C)} $ has a minimal generator of degree $ 3 $ if and only if $C$ is either trigonal or isomorphic to a plane quintic.

On the other hand if $\im(\psi_C)$ is the orbit of the Fermat cubic, then  $\im(\psi_C)\subseteq H_{g-3,3}(g-2)$. Since Lemma 1.5 implies $\im(\psi_C)\not\subseteq H_{g-3,3}(g-3)$, it follows that $\ap(C)=g-2$. Conversely let us assume that $\ap(C)=g-2$, so that the general point of $\im(\psi_C)$ is in $H_{g-3,3}(g-2)\setminus H_{g-3,3}(g-3)$. If the general point of $\im(\psi_C)$ is in the orbit of $f:=\sum_{i=0}^{g-3}\lambda_i\ell_i^{[3]}$, then the linear forms $\ell_0,\dots,\ell_{g-3}$ are linearly independent by Lemma 1.5, thus $f$ is a, possibly divided power, Fermat cubic up to $\GL_{g-2}$.
\qed
\enddemo

\remark{Remark 1.8}
In the proof of Theorem 1.7 above we saw that
$$
f^{\perp} =(x_i x_j , x_i^3-x_{g-3}^3)_{0 \le i < j \le g-3}.
$$
Let $\Gamma\subseteq H:=D_0(x_{g-2},x_{g-1})$ consist of the $g-2$ fundamental points
$$
E_0:=[1,0,\dots,0,0,0,0],\quad\dots,\quad E_{g-3}:=[0,0,\dots,0,1,0,0]
$$
in $H\subseteq\dualp{g-1}$. In what follows we will show that $ T^{g-3}/ f^{\perp} $ can be obtained from $ T({\Gamma})$ via the anticanonical construction (see [Mi], Theorem 4.2.8).

The homogeneous ideal of $\Gamma$ in $T^{g-3}$ is $I_{\Gamma\vert H}=(x_i x_j)_{0 \le i < j \le g-3}$,
and the Betti numbers of its minimal free resolution are the same as the defining ideal of the rational normal curve in $ T^{g-2}$.

In order to describe the maps of such a resolution, for each $h=1,\dots,g-3$ let us consider the natural exact sequence 
$$
0 \longrightarrow I_h\cap (x_{h+1})J_h\longrightarrow I_h \oplus (x_{h+1})J_h \longrightarrow I_{h}+(x_{h+1}) J_h\longrightarrow 0 
$$ 
where $ I_h :=(x_i x_j)_{0 \le i < j \le h}$, $ J_h :=(x_0, \dots, x_h)$, the first non--trivial map is $u\mapsto(u,-u)$, the second one is $(u,v)\mapsto u+v$. We trivially have $ I_h(-1)\cong x_{h+1} I_h = I_h\cap (x_{h+1} )J_h$, $ I_h+ (x_{h+1})J_h = I_{h+1}$ and $ (x_{h+1}) J_h \cong J_h(-1)$ (here and in the following $(-1)$ denotes the degree shifting, as usual), then the above sequence is isomorphic to a sequence of the form 
$$
0 \longrightarrow I_h(-1) \longrightarrow I_h \oplus J_h(-1) \longrightarrow I_{h+1} \longrightarrow 0 
$$ 
where the first non--trivial map is $u\mapsto(x_{h+1}u,-u)$ and the second one is $(u,v)\mapsto (u+x_{h+1}v)$.

The last above sequence allows us to compute inductively via mapping cone the maps of the minimal free resolution of $I_{\Gamma\vert H}$. In particular such a resolution ends with
$$
0\longrightarrow(T^{g-3})^{\oplus (g-3)}(-g+2) \mapright{\delta} (T^{g-3})^{\oplus (g-4)(g-2)}(-g+3).
$$
In order to describe the matrix $M:=(m_{p,q})$ of $\delta$ we set $ X_i = (x_{g-3}, \dots, x_{i+1}, x_{i-1}, \dots, x_0) $, $ i = 0, \dots, g-5$, and $ X_{g-4} = (x_{g-5}, \dots, x_0)$. Then we have
$$ 
m_{p,q} = \cases
(-1)^{q+j} (X_{q-1})_j & \text{if $p = (q-1)(g-3) + j$, $1 \le q \le g-4$, $1 \le j \le g-3$}, \\
 (-1)^{g-3+j} (X_{g-4})_j & \text{if $p = (g-4)(g-3) + j$, $q = g-3$, $1 \le j \le g-4$} ,\\ 
 (-1)^{i+1} x_{g-3} & \text{if $p = i(g-3) + 1$, $q = g-3$, $0 \le i \le g-5$}, \\ 
 0 & \text{othewise} ,
 \endcases
$$
where $ (X_s)_j $ is the $j$-th entry of the vector $ X_s$.

Let $ K_{\Gamma} $ be the canonical module of $ \Gamma$, that is to say, $ K_{\Gamma} = Ext_S^{g-3} \big(T({{\Gamma}}), T^{g-3} \big) (-g+2)$. The scheme $\Gamma$ is aCM since it is $0$--dimensional, hence the minimal free resolution of $ K_{\Gamma} $ as module over $T^{g-3}$ can be obtained by dualizing the one of $ I_{\Gamma\vert H}$. In particular, the minimal free resolution of $ K_{\Gamma} $ starts as
$$ 
(T^{g-3})^{\oplus (g-4)(g-2)} (-1) \mapright{\check\delta}(T^{g-3})^{\oplus (g-3)} \longrightarrow K_{\Gamma}\longrightarrow0,
$$
where $\check\delta$ is represented by $ {}^tM$, the transpose of $ M$.

We want to prove that $ f^{\perp}/I_{\Gamma\vert H} \cong K_{\Gamma}(-3)$. At first, we have the obvious short exact sequence 
$$
0 \longrightarrow {f^{\perp}}/{I_{\Gamma\vert H}}  \longrightarrow T({{\Gamma}}) \longrightarrow T^{g-3}/{f^{\perp}}  \longrightarrow 0.
$$
By comparing the minimal free resolutions of the last two modules, and the fact that $ f^{\perp}/I_{\Gamma\vert H} $ is generated in degree $ 3$, we get that the Betti numbers of $ K_{\Gamma}(-3) $ and $ f^{\perp}/I_{\Gamma\vert H} $ are equal. Moreover one can also check directly that the columns of $ {}^tM$ are syzygies of $ x_0^3 - x_{g-3}^3, \dots, x_{g-4}^3 - x_{g-3}^3 $. Hence, both $ K_{\Gamma}(-3) $ and $ f^{\perp}/I_{\Gamma\vert H} $ are minimally presented by $ {}^tM$ and so the isomorphism follows.

Thus we have finally proved that $ T^{g-3}/ f^{\perp} $ is obtained from $ T({\Gamma})$ via the anticanonical construction (see [Mi], Theorem 4.2.8) as stated at the begining of this remark.
\endremark
\medbreak

\head
2. Surfaces of almost minimal degree
\endhead
In this section we will recall some facts about canonical curves of genus $g$ lying on integral surfaces of degree $g-1$ in the canonical space, i.e. surfaces of almost minimal degree (see the introduction). Such results are more or less known  (see e.g. [Sch2], Section 4).

We start by recalling the definition and the classification of surfaces of almost minimal degree.

\proclaim{Theorem 2.1}
Let $S\subseteq\p n$ be an integral non--degenerate surface of degree $n\ge4$. Then the following assertions hold true.
\item{i)} If $S$ is not aCM or it is aCM but not normal, then it is the projection of a  surface of minimal degree in $\p{n+1}$ via a linear projection.
\item{ii)} If $S$ is both aCM and normal, then $\omega_S\cong\Ofa S(-1)$: in particular $S$ is aG. 
\endproclaim
\demo{Proof}
See [B--S], Theorem 6.2 and Corollary 6.10. In [B--S] a very explicit description of the structure of the projection when $S$ is not aCM or it is aCM but not normal can be also found. 
\qed
\enddemo

In the next section we will see that the class of integral normal (hence aG and anticanonically embeded due to the above Theorem) surfaces is actually the most interesting for us. For this reason we now go to summarize some facts about them (see [H--W] for their proofs).

\remark{Remark 2.2}
If $S$ is an integral normal cone, then the general hyperplane section $E$ of $S$ is an aCM curve of degree $n$ in $\p{n-1}$. Thus Castelnuovo's bound (see [A--C--G--H], Chapter III, Section 2 when $\char(k)=0$ and [Ha] Theorem IV.6.4 or [Ra], Theorem 2.9 when $\char(k)>0$) yields that the arithmetic genus of $E$ is either $0$ or $1$. Since $E$ is aCM, it follows that it must be necessarily an elliptic curve.
\endremark
\medbreak

Now let us assume that $S$ is not a cone. In this case in [H--W] it is proved that such integral surfaces coincide with the ones described in [De]. Before stating the main theorem, we recall some facts about them. Following the notation and definitions introduced in [H--W] and [De], let $r\le8$ be a positive integer and consider (possibly infinitely near) points  $p_1,\dots,p_r$ such that: $p_1\in \widetilde{S}_0:=\p2$, $p_{j+1}\in \widetilde{S}_j$, $j=1,\dots,r-1$ and $\widetilde{S}_j$ is the blow up of $\widetilde{S}_{j-1}$ at $p_{j}$, $j=1,\dots,r$.
We thus have a chain of blow up's
$$
\widetilde{S}:=\widetilde{S}_r\to \widetilde{S}_{r-1}\to\dots\to \widetilde{S}_2\to \widetilde{S}_1\to \widetilde{S}_0=\p2.
$$
We will denote by $E_j\subseteq\widetilde{S}_j$ the exceptional divisor of the  blow up $\widetilde{S}_j\to \widetilde{S}_{j-1}$.

\definition{Definition 2.3}
With the above notation the set $\{\ p_1,\dots,p_r\ \}$ is said to be in almost general position if the following conditions hold.
\item{i)} No four of them are on the same line.
\item{ii)} No seven  of them are on the same conic.
\item{iii)} For each $j=1,\dots,r-1$ the point $p_{j+1}\in \widetilde{S}_j$  does not lie on any proper transform $\widehat{E_i}$ of $E_i$, $i=1,\dots,j$, such that $\widehat{E_i}^2=-2$.
\enddefinition

Let $\widetilde{S}$ be a surface obtained by blowing up a set of points in almost general position. The Picard group $\Pic(\widetilde{S})$ of $\widetilde{S}$ is freely generated by the class $\ell$ of the proper transform of a line in $\p2$ and by the $r$ classes $e_1,\dots,e_r$ of the total transforms of the exceptional divisors $E_1,\dots,E_r$. Notice that $\ell^2=1$, $\ell\cdot e_i=0$, $e_i^2=-1$ and $e_i\cdot e_j=0$, $i,j=1,\dots,r$, $i\ne j$. The canonical system on $\widetilde{S}$ is $\vert3\ell-\sum_{i=1}^re_i\vert$. The anticanonical map is a birational morphism onto a surface $S\subseteq\p{9-r}$ of degree $9-r$ when $r\le6$. In [H-W] a converse of this fact is proven. We write below the part of this result that we will need in the following.

\proclaim{Theorem 2.4}
Let $S\subseteq\p n$ be an integral normal and linearly normal surface of degree $n $. Assume that $S$ is not a cone on a curve of degree $n$ lying in a hyperplane. Then the following assertions hold true.
\item{i)} $3\le n\le 9$.
\item{ii)} $S$ carries at most rational double points as singularities.
\item{iii)} If $n=9$, then $S$ is the anticanonical image in $\p9$ of $\p2$.
\item{iv)} If $n=8$ then $S$ is the anticanonical image in $\p8$ of either $\p1\times\p1$ or ${\Bbb F}_1:={\Bbb P}(\Ofa{\p1}\oplus\Ofa{\p1}(-1))$ or ${\Bbb F}_2:={\Bbb P}(\Ofa{\p1}\oplus\Ofa{\p1}(-2))$.
\item{v)} If $3\le n\le7$, then $S$ is the anticanonical image in $\p n$ of the blow up of $\p2$ along a set of $9-n$ points in almost general position.

\noindent In all the above  cases the blow up map $\widetilde{S}\to S$ is the contraction of all the curves $D\subseteq \widetilde{S}$ such that $D^2=-2$, if any.
\endproclaim
\demo{Proof}
See [H--W], Theorem 3.4, for the general result.
\qed
\enddemo

\definition{Definition 2.5}
The embedded surfaces described in the previous statement will be called weak del Pezzo surfaces of degree $n$.
\enddefinition

\head
3. Canonical curves on surfaces of almost minimal degree
\endhead
We recall here some facts about the canonical embedding of non--hyperelliptic curves (recall that twe are considering only smooth curves). First of all, let  $C$ be a non--hyperelliptic curve of genus $g\ge5$. Then $I_{\phi(C)}$ has a minimal free resolution over $T^{g-1}$ of the form
$$
0\longrightarrow F_{g-2}\longrightarrow F_{g-1}\longrightarrow\dots\longrightarrow
F_1\longrightarrow I_{\phi(C)}\longrightarrow0.\tag 3.1
$$
Since the embedding $\phi(C)\subseteq\dualp{g-1}$ is aG, such a resolution is isomorphic to its dual up to twist. Thus 
$$
F_p\cong T^{g-1}(-p-1)^{\oplus\beta_{p,p+1}}\oplus T^{g-1}(-p-2)^{\oplus\beta_{p,p+2}}
$$
if $p=1,\dots,g-3$ and $F_{g-2}\cong T^{g-1}(-g-1)$. In particular the Hilbert polynomial of $T_C$ is
$$
\align
{{g-1+t\choose g-1}}&+\sum_{p=1}^{g-3}(-1)^p\left(\beta_{p,p+1}{{g-2+t-p\choose g-1}}+\right.\\
&\left.+\beta_{p,p+2}{{g-3+t-p\choose g-1}}\right)+(-1)^{g-2}{{t-2\choose g-1}}.
\endalign
$$
It follows that
$$
\beta_{p,p+1}-\beta_{p-1,p+1}=p{g-2\choose{p+1}}-(g-1-p){g-2\choose{g-p}}.\tag 3.2
$$
We conclude this preliminary part with the following

\proclaim{Lemma 3.3}
Let $C$ be a non--hyperelliptic curve of genus $g\ge5$. Assume that the canonical model $\phi(C)$ of $C$ is contained in an integral aCM normal surface $S\subseteq\dualp{g-1}$ of degree  $g-1$. Then $\phi(C)$ is cut out on $S$ by a quadric hypersurface $Q$ such that $\Sing(S)\cap Q=\emptyset$. In particular $C$ is neither trigonal nor isomorphic to a plane quintic.
\endproclaim
\demo{Proof}
By Theorem 2.1 the surface $S$ is aG. Up to a proper choice of the coordinates $x_0,\dots,x_{g-1}$ in $\dualp{g-1}$ we can assume that  $(x_{g-2},x_{g-1})$ is a regular sequence in $T(S)$. Thus the subspace $H:=D_0(x_{g-2},x_{g-1})\subseteq\dualp{g-1}$ has codimension $2$ and the minimal free resolution of the homogeneous coordinate ring $T(\Gamma)$ of the scheme $\Gamma:=S\cap H$ over $T^{g-1}/(x_{g-2},x_{g-1})\cong T^{g-3}$ has the same Betti numbers of the minimal free resolution of the $T^{g-1}$--module $T(S)$ (the argument is the same used in Remark 1.3). 

In [G--O] (see also Section 4 of [Sch1]) the Betti numbers of  $I_{\Gamma\vert H}$ are computed. In particular we know that its minimal free resolution of $ I_{\Gamma\vert H}$ has the shape
$$
0\longrightarrow T^{g-3}(1-g)\longrightarrow T^{g-3}(3-g)^{\oplus\gamma_{g-4}}\longrightarrow\dots\longrightarrow
T^{g-3}(-2)^{\oplus\gamma_1}\longrightarrow I_{\Gamma\vert H}\longrightarrow0\tag3.3.1
$$
where
$$
\gamma_i:={{i(g-3-i)}\over{g-2}}{g-1\choose{i+1}}.
$$
It follows from the above sequence and from Formula (3.2) above that $I_S$ is generated by 
$$
\gamma_1:={g-2\choose 2}-1=\beta_{1,2}-1.
$$
In particular, there exists a quadric $Q\subseteq\dualp{g-1}$ such that $\phi(C)\subseteq S\cap Q$. Since $S$ is integral, due to degree reasons, equality must hold. The tangent space $T_p(C)$ of $C$ at any point $p\in C$ is the intersection $T_p(S)\cap T_p(Q)$ of the tangent spaces of $S$ and $Q$ at the same point $p$. Since $Q$ is a hypersurface it follows that necessarily $\dim(T_p(S))=2$ for each $p\in C$.

It is now an easy consequence of the Babbage--Enriques--Noether--Petri Theorem (see [SD]) that $C$ is neither trigonal nor isomorphic to a plane quintic.
\qed
\enddemo

We are now ready to recall the main result of this section which is a quite natural generalization of the classical result for canonical curves on surfaces of minimal degree. 

\proclaim{Theorem 3.4}
Let $C$ be a non--hyperelliptic  curve of genus $g\ge5$.
Then the canonical model $\phi(C)$ of $C$ is contained in an integral surface $S\subseteq\dualp{g-1}$ of degree  $g-1$ if and only if one of the three following conditions holds.
\item{i)} $C$ is bielliptic: in this case $S$ is a cone on an elliptic normal curve contained in a hyperplane of $\dualp{g-1}$ with vertex not on $\phi(C)$.
\item{ii)} $g\le10$ and $C$ is birationally isomorphic to a plane sextic carrying $10-g$ double points in almost general position as singularities: in this case $S$ is the weak del Pezzo surface obtained by embedding anticanonically $\p2$ blown up at the singular points of the above mentioned plane sextic.
\item{iii)} $g=9$ and $C$ is isomorphic to a smooth complete intersection inside $\p3$ of an integral quadric $Q$ with an integral quartic surface: in this case $S$ is the weak del Pezzo surface obtained by embedding anticanonically $Q$.
\endproclaim
\demo{Proof}
We first prove the \lq\lq if\rq\rq\ part of the statement. If $C$ is birationally isomorphic to a plane sextic carrying $10-g$ double points in almost general position as singularities, then $C$ lies on the weak del Pezzo surface $S$ obtained as anticanonical embedding of the blow up of $\p2$ at the singular points of the above mentioned plane sextic. A similar argument holds if $C$ is isomorphic to a smooth complete intersection inside $\p3$ of an integral quadric $Q\subseteq\p3$ with an integral quartic surface. 

Now we examine the case when $C$ is bielliptic. For the following construction we refer to [Sch1] (in particular see Section 6). Let $C$ be bielliptic. Thus there exists a morphism $\varphi\colon C\to E$ of degree $2$ onto an elliptic curve. The pull--back to $C$ of each $g^1_2$ on $E$ gives rise to a $g^1_4$ on $C$. Fix one of such $g^1_4$'s: the union of planes in the canonical space $\dualp{g-1}$ generated by the divisors belonging to the fixed $g^1_4$  is a singular rational normal scroll $\Sigma\subseteq\dualp{g-1}$ over $\p1$. 

The canonical model $\phi(C)$ of $C$ is the complete intersection inside $\Sigma$ of two surfaces, namely a cone $S$ of degree $g-1$ and another surface of degree $2g-6$. Since the curve $\phi(C)$ is smooth, it follows that the two surfaces are necessarily smooth at the points of $\phi(C)$. Let us consider the general hyperplane section $E$ of $S$. The scheme $E$  has dimension $1$, degree $g-1$ and it lies in a hyperplane of $\dualp{g-1}$. Moreover it is non--degenerate since, otherwise, the curve $\phi(C)$ would be degenerate. If $E$ were singular, then $S$, being a cone on $E$,  would have a line $r$ of singular points, thus $S$ would be singular at the points in $r\cap\phi(C)$, a contradiction.

Thus Castelnuovo's bound (see [A--C--G--H], Chapter III, Section 2 when $\char(k)=0$ and [Ha] Theorem IV.6.4 or [Ra], Theorem 2.9 when $\char(k)>0$) yields that the (geometric) genus of $E$ is either $0$ or $1$ and, in the first case $C$ would be hyperelliptic. We conclude that $E$ is an elliptic curve. Due to Corollary 2.3 of [C--G--N] it follows that $E$ is aCM, thus the same is true for $S$.

We have thus proved that in all the above cases the canonical model of $C$ lies on an integral surface $S\subseteq\dualp{g-1}$ of degree $g-1$. Moreover we have also given a complete description of such a surface $S$.

Now we prove the \lq\lq only if\rq\rq\ part of the statement.
Let us assume that $\phi(C)$ is contained in an integral surface $S\subseteq\dualp{g-1}$ of degree  $g-1$. Since $\phi(C)\subseteq\dualp{g-1}$ is non--degenerate, the same holds for $S\subseteq\dualp{g-1}$.

If $S$ were not aCM or it is aCM but not normal, then $S$ would be obtained as the linear projection of a surface $\Sigma$ of degree $g-1$ contained in a projective space  of dimension $g\ge6$ and containing $\dualp{g-1}$ as a hyperplane. In [E--H], Theorem 1, it is proved that such a surface $\Sigma$ would either be a ruled (possibly singular) surface or the Veronese surface in $\p5$ (hence $g=5$).

Consider the first case. The images of the fibres of $\Sigma$ are lines sweeping $S$ and cutting out a $g^1_d$ on $\phi(C)$. 
Necessarily $d\ge3$ since $C$ is neither rational nor hyperelliptic. Let $\vert D\vert$ be the complete linear system containing such a $g^1_d$. The geometric interpretation of Riemann--Roch Theorem (see [A--C--G--H], Chapter 1, Section 1: the formula holds in any characteristic) implies that the spaces generated by such divisors in $\dualp{g-1}$ have dimension $d-1-\dim\vert D\vert$, whence $\dim\vert D\vert=d-2$. Since $g\ge5$ such a linear system is special, thus Clifford's theorem yields
$d-2\le d/2$, i.e. $d\le4$. If equality holds, then $\vert D\vert$ would be the canonical system, thus  $d=3$, whence $C$ is trigonal. In Theorem 4.7 and Lemma 4.8 of [SD] it is proved that the quadrics through $\phi(C)$ cut out a surface $S'\subseteq\dualp{g-1}$ of degree $g-2$. Since $S'$ is intersection of quadrics, we deduce that it contains each line cutting the $g^1_3$ on $\phi(C)$, hence $S\subseteq S'$, a contradiction. 

In the second case $\phi(C)$ would be the birational image of an integral divisor $F\subseteq V\subseteq\p5$ via a projection from a point not on $V$, hence $\deg(F)=\deg(\phi(C))=8$ and its genus should be $5$. Since $V\cong\p2$, it follows that $F\in\vert\Ofa{\p2}(4)\vert$, hence the arithmetic genus of $F$ would be at most $3$, a contradiction.

Thus we can assume that $C$ is both aCM and normal. Due to the results listed in Section 2 (in particular see Remark 2.2) the surface $S$ is either a weak del Pezzo surface or it is a cone over an elliptic curve $E$ in a hyperplane $H\subseteq\dualp{g-1}$ with vertex a point $V\not\in H$. 

In both the cases $\phi(C)$ is cut out on $S$ by a single quadric (see Lemma 3.3) not containing any singular point of $S$. In particular $C\in\vert\omega_S^{-2}\vert$. Thus in the first case, thanks to the classification of weak del Pezzo surfaces given in Theorem 2.4, we obtain curves which are either birationally isomorphic to plane sextics carrying $10-g$ double points in almost general position as singularities or that are isomorphic to smooth complete intersections inside $\p3$ of integral quadrics $Q$ with an integral quartics. 
In the second case the vertex $V$ of the cone $S$ is not on the quadric $Q$, thus the projection $\phi(C)\to E$ with center $V$ has degree $2$, i.e. $C$ is bielliptic.
\qed
\enddemo

\head
3. Curves with almost minimal apolarity
\endhead
In this section we will prove Theorem B stated in the introduction. To this purpose we first prove the following technical lemma.

\proclaim{Lemma 4.1}
Let $C$ be a non--hyperelliptic curve of genus $g\ge5$ apolar to the polynomial $f\in R^{g-3}_3$. If $f=f_0+f_1y_{g-3}^{[3]}$ where $f_0\in R^{g-4}_3\subseteq R^{g-3}_3$ and $f_1\in k$, then $C$ is either trigonal or isomorphic to a plane quintic.
\endproclaim
\demo{Proof}
In order to show that $C$ is trigonal or isomorphic to a plane quintic it suffices to prove that a minimal set of generators of the homogeneous ideal $I_{\phi(C)}\subseteq T^{g-1}$ of its canonical model $\phi(C)\subseteq\dualp{g-1}$ contains at least a non--quadratic polynomial. Since the Betti numbers of $I_{\phi(C)}\subseteq T^{g-1}$ coincide with the ones of $f^\perp\subseteq T^{g-3}$ (see Remark 1.3), it then suffices to prove that a minimal set of generators of $f^\perp$ contains at least a non--quadratic polynomial.

Choose a quadratic polynomial $q\in f^\perp$: we can write
$$
q=q_0+\sum_{i=0}^{g-4}a_ix_ix_{g-3}+a_{g-3}x_{g-3}^2
$$
where $q_0\in R^{g-4}_2$ and $a_0,\dots,a_{g-3}\in k$. Since we certainly have $f_1\ne0$ (Lemma 1.5), the equality
$$
 q(f)=q_0(f_0)+a_{g-3}f_1y_{g-3}
$$
yields $a_{g-3}=0$. We conclude that each quadratic polynomial $q\in f^\perp$ vanishes at $[0,\dots,0,1]$. Since $T^{g-3}/f^\perp$ has dimension $0$ it follows that there is another minimal non--quadratic generator of $f^\perp$.
\qed
\enddemo

We are now ready to prove Theorem B stated in the introduction.

\proclaim{Theorem 4.2}
A non--hyperelliptic curve $C$ of genus $g\ge5$ is either bielliptic or birationally isomorphic to a plane sextic carrying $10-g$ double points in almost general position as singularities or isomorphic to a smooth complete intersection inside $\p3$ of an integral quadric $Q$ with an integral quartic surface if and only if $\ap(C)=g-1$.
\endproclaim
\demo{Proof}
In view of Theorem 3.4 it suffices to prove that $\ap(C)=g-1$ if and only if there exists an integral surface $S\subseteq\dualp{g-1}$ of degree $g-1$ containing the canonical model $\phi(C)$ of $C$. 

Assume that such a surface $S\subseteq\dualp{g-1}$ exists. Thus $S$ is either a weak del Pezzo surface or it is a cone over an elliptic normal curve contained in a hyperplane. In both the cases $S$ has at most a finite number of singular points due to Theorem 3.4.

Let $H:=D_0(h_1,h_2)\subseteq\dualp{g-1}$ where  $h_1$ and $h_2$ are general linear forms. The geometric characterization of the degree of a surface allows us to assume that $\Gamma:=S\cap H$ consists of exactly $g-1$ pairwise distinct points. 

As usual let $I_{S}$ be the homogeneous ideal of $S$. We know that $I_\Gamma$ is the saturation of $I_{S}+(h_1,h_2)$. By construction we have 
$$
{{I_{S} +(h_1,h_2)}\over{(h_1,h_2)}}\subseteq{{I_{\phi(C)} +(h_1,h_2)}\over{(h_1,h_2)}}
$$
and the quotient on the right is saturated, since the same is true for $I_{\phi(C)}$, being $\phi(C)$ aG (see [Mi], Proposition 1.3.4). Thus, by saturating both the sides of the above inclusion the quotient on the right does not change and we finally obtain 
$$
I_{\Gamma\vert H}\subseteq {{I_{\phi(C)} +(h_1,h_2)}\over{(h_1,h_2)}}.
$$
We can assume that $h_i=x _{g-3+i}$, $i=1,2$. Via the natural identification $T^{g-1}/(h_1,h_2)=T^{g-3}$, we can find $f\in R^{g-3}_3$ such that
$$
f^\perp= {{I_{\phi(C)} +(h_1,h_2)}\over{(h_1,h_2)}}
$$
(as in Section 1), hence $I_{\Gamma\vert H}\subseteq f^\perp\subseteq {T^{g-3}}$ and a natural identification
$$
{{T_C}\over{(h_1,h_2)}}={{T^{g-3}}\over{f^\perp}}.
$$
It then follows from Lemma 1.6 that the general point of $\im(\psi_C)$ is in $H_{g-3,3}(g-1)$, hence $\ap(C)\le g-1$. If $\ap(C)<g-1$, then $C$ would be either trigonal or isomorphic to a plane quintic due to Theorem 1.7, a contradiction due to Lemma 3.3.

Let us now prove the converse. To this purpose let us assume that $\ap(C)=g-1$: hence $C$ is neither trigonal nor isomorphic to a plane quintic. Furthermore we can find linear forms $\ell_0,\dots,\ell_{g-2}\in T^{g-3}$ such that 
$f=\sum_{i=0}^{g-2}\ell_i^{[3]}$. 

Due to Lemma 1.5 we can assume that $\ell_i=y_i$ and $\ell_{g-2}=\sum_{i=0}^{g-3}\lambda_iy_i$, where $\lambda_i\in k$, $i=0,\dots,g-3$. Lemma 4.1 forces the non--vanishing of $\lambda_i$, $i=0,\dots,g-3$. The transformation $y_i\mapsto y_i/\lambda_i$, $i=0,\dots,g-3$ allows us to assume 
$$
f:=\sum_{i=0}^{g-3}\lambda_i^{-1}y_i^{[3]}+\left(\sum_{i=0}^{g-3}y_i\right)^{[3]}.\tag4.2.1
$$
It follows by the Apolarity Lemma that $I_{\Gamma\vert H}\subseteq f^\perp$, where $I_{\Gamma\vert H}$ is the ideal of  the set $\Gamma\subseteq H$ consisting of the $g-1$ fundamental points
$$
E_0:=[1,0,\dots,0,0,0,0],\quad\dots,\quad E_{g-3}:=[0,0,\dots,0,1,0,0],\quad E_{g-2}:=[1,1,\dots,1,1,0,0]
$$
giving the standard projective frame in $H\subseteq\dualp{g-1}$. 

We thus have that the minimal free resolution of $I_{\Gamma\vert H}$ over $T^{g-3}$ coincides with Resolution (3.3.1), while the Betti numbers of the minimal free resolution of $f^\perp$ over $T^{g-3}$ coincide with those in Resolution (3.1) due to Remark 1.3. The fact that $C$ is neither trigonal nor isomorphic to a plane quintic yields that $\beta_{1,3}=0$, hence Formula (3.2) with $p=1$ implies that $f^\perp$ is generated by
$$
\beta_{1,2}={g-2\choose{2}}=\gamma_{1}+1
$$
quadratic polynomials. It follows that $f^\perp=I_{\Gamma\vert H}+(q_f)$ for a suitable $q_f\in T^{g-3}_2$. Fix an ordered minimal set of generators $(q_1,\dots,q_{\gamma_{1}},q_f)$ of $f^\perp$ in such a way that $I_{\Gamma\vert H}=(q_1,\dots,q_{\gamma_{1}})$

Moreover Formula (3.2) with $p=2$ gives
$$
\beta_{2,3}={{(g-1)(g-3)(g-5)}\over3}=\gamma_{2},
$$ 
hence the linear syzygies of $f^\perp$ are exactly the linear syzygies of $I_{\Gamma\vert H}$ with a $0$ in correspondence of the last generator $q_f$. In degree $2$ we have the $\gamma_{1}={g-2\choose{2}}-1$
Koszul syzygies of the form
$$
(0,\dots, 0,-q_f,0,\dots,0,q_i)
$$
where the $q_f$ sits in the $i^{\roman th}$ position. Such syzygies are obviously not generated by the syzygies in degree $1$ since their last entries are always non--zero, thus $\beta_{2,4}=\gamma_{1}$. It follows now from the table in Theorem 4.1 of [Sch2] that $\phi(C)$ lies on an integral surface of degree $g-1$. Since $C$ is neither trigonal nor isomorphic to a plane quintic then there are no surfaces of lower degree containing $\phi(C)$.
\qed
\enddemo

\remark{Remark 4.3}
We know by the proof of Theorem 4.2 that we can assume that $\Gamma=S\cap H$ consists of the $g-1$ fundamental points $E_0,\dots,E_{g-3}, E_{g-2}$ giving the standard projective frame, thus 
$$
I_{\Gamma\vert H}=\left(x_hx_i-x_hx_j\right)_{0\le h<i< j\le g-3}.
$$
Moreover we also know that $f^\perp=I_{\Gamma\vert H}+(q_f)$, thus it could be interesting to identify the class of $q_f$ modulo $I_{\Gamma\vert H}$. To this purpose we can assume that 
$$
f:=\sum_{i=0}^{g-3}\lambda_i^{-1}y_i^{[3]}+\left(\sum_{i=0}^{g-3}y_i\right)^{[3]}
$$
for suitable non--zero $\lambda_i$, $i=0,\dots,g-3$ (see Equality (4.2.1)). It is easy to check that the polynomial 
$$
q':=x_{g-4}x_{g-3}-\sum_{i=0}^{g-3}{\lambda_i}\left(x_i^2-x_{g-4}x_{g-3}\right)
$$
is in $f^\perp$. Since $q'$ does not vanish at $E_{0}=[1,0,\dots,0]$, it follows that $q'\not\in I_{\Gamma\vert H}$, thus $q_f\equiv q'$ modulo $I_{\Gamma\vert H}$.
\endremark
\medbreak

\remark{Remark 4.4}
The quotient morphism
$$
T^{g-1}\twoheadrightarrow T^{g-1}/(h_1,h_2)\cong T^{g-3}
$$
induces by restriction an isomorphism $\psi\colon(I_{\phi(C)})_2\mapright\sim (f^\perp)_2$ of $k$--vector spaces. Recall that the subscheme $\Gamma\subseteq H$ satisfies $I_{\Gamma\vert H}\subseteq f^\perp$. Thus the ideal $I_S$ of $S$ is the ideal generated by $\psi^{-1}((I_{\Gamma\vert H})_2)$.
\endremark
\medbreak

\remark{Remark 4.5}
With the same notation used in the proof of Theorem 4.2 and in Remarks 4.3, 4.4, equality $ I_{\Gamma \vert H} \cap (q_f ) = I_{\Gamma \vert H} \cdot (q_f )$ holds. In fact, if $ q_f g \in I_{\Gamma \vert H}$, then $ q_f(E_i) g(E_i) = 0 $ for each point $ E_i$, with $ i = 0, \dots, g-2$. From the equality $ I_{\Gamma \vert H} + (q_f) = f^{\perp} $ and the fact that $ T^{g-3} / f^{\perp} $ has Krull--dimension $ 0$, we deduce that $ g(E_i) = 0 $ for each $ i = 0, \dots, g-2 $ i.e. $ g \in I_{\Gamma \vert H} $ and the claim follows.

From the previous equality, we deduce the exactness of the following short sequence of modules over $T^{g-1}/(x_{g-2},x_{g-1})\cong T^{g-3}$
$$
0 \longrightarrow T(\Gamma)(-2) \mapright{\cdot q_f} T(\Gamma)\longrightarrow T^{g-3}/{f^{\perp}}  \longrightarrow 0 .
$$ 
It shows that we can get the minimal free resolution of $ T^{g-3}/{f^{\perp}}  $ via mapping cone from the one of $ T(\Gamma)$.

Furthermore, let $ K_{\Gamma} $ be the canonical module of $ \Gamma$, that is to say, 
$$
K_{\Gamma} := \Ext^{g-3}_{T^{g-3}} \big(T(\Gamma), T^{g-3}\big)(-g+2).
$$
From the minimal free resolution of $ I_{\Gamma \vert H} $ (see Resolution (3.3.1)) and standard results on Gorenstein rings, we get that $ K_{\Gamma}(-1) \cong T(\Gamma)$. From such an equality, it then follows that 
$$
{f^{\perp}}/{I_{\Gamma \vert H}} \cong T(\Gamma)(-2) \cong K_{\Gamma}(-3),
$$
that is to say, $ T^{g-3}/{f^{\perp}} $ can be obtained from $T(\Gamma)$ via the  anticanonical divisor construction (see [Mi], Theorem 4.2.8).
\endremark
\medbreak


\Refs
\refstyle{A}
\widestnumber\key{A--C--G--H}

\ref
\key A--H
\by J\. Alexander, A\. Hirschowitz
\paper Polynomial interpolation in several variables
\jour J. Algebraic Geom. 
\vol 4 
\yr 1995
\pages  201--222
\endref

\ref
\key A--C--G--H
\by E\. Arbarello, M\. Cornalba, P.A\. Griffiths, J\. Harris
\book Geometry of algebraic curves
\vol I
\publ Springer
\yr 1985
\endref


\ref
\key B--S
\by M\. Brodmann, P\. Schenzel
\paper Arithmetic properties of projective varieties of almost minimal degree
\jour J. Algebraic Geom.
\vol 16
\yr 2007
\pages 347--400
\endref

\ref
\key C--N1
\by G\. Casnati, R\. Notari
\paper On the Gorenstein locus of some punctual Hilbert schemes
\jour J. Pure Appl. Algebra
\vol 213
\yr 2009
\pages 2055-2074
\endref

\ref
\key C--N2
\bysame % G\. Casnati, R\. Notari
\paper On the irreducibility of the Gorenstein locus of the punctual Hilbert schemes of degree $10$
\paperinfo in preparation
\endref

\ref
\key C--G--N
\by N\. Chiarli, S\. Greco, U\. Nagel
\paper On the genus and Hartshorne--Rao module of projective curves
\jour Math. Z.
\vol 229
\yr 1998
\pages  695--724
\endref

\ref
\key C--R--V
\by A\. Conca, M.E\. Rossi, G\.Valla
\paper Gr\"obner flags and Gorenstein algebras
\jour Compositio Math.
\vol 129
\yr 2001
\pages 95--121
\endref

\ref
\key DP--Z1
\by  P\. De Poi, F\. Zucconi
\paper Gonality, apolarity, and hypercubics
\paperinfo arXiv: 0802.0705v2 [math.AG]
\endref

\ref
\key DP--Z2
\bysame %P\. De Poi, F\. Zucconi
\paper Fermat hypersurfaces and Subcanonical curves
\paperinfo arXiv: 0908.0522v1  [math.AG]
\endref

\ref
\key De
\by M\. Demazure
\paper Surfaces de Del Pezzo
\inbook S\'eminaire sur les singularit\'es des surfaces (Palaiseau, France 1976--1977)
\eds M. Demazure, H. Pinkham, B. Teissier
\bookinfo L.N.M.
\vol 777
\yr 1980
\pages 23--69
\endref


\ref
\key E--H
\by D\. Eisenbud, J\. Harris
\paper On varieties of minimal degree (a centennial account)
\inbook Algebraic geometry, Bowdoin 1985
\ed Spencer J.  Bloch
\bookinfo Proceedings of symposia in pure mathematics
\vol 46, part I
\pages 3--13
\publ A.M.S.
\yr 1987
\endref

\ref
\key G--O
\by A.V\. Geramita, F\. Orecchia
\paper On the Cohen--Macaulay type of $s$ lines in ${\Bbb A}^{n+1}$
\jour J. Alg.
\vol 70
\yr 1981
\pages  116--140
\endref

\ref
\key Ha
\by R\. Hartshorne
\book Algebraic geometry
\bookinfo G.T.M. 
\vol 52
\publ Springer
\yr 1977
\endref

\ref
\key H--W
\by F\. Hidaka, K\. Watanabe
\paper Normal Gorenstein surfaces with ample anti-canonical divisor
\jour Tokyo J. Math.
\vol 4
\yr 1981
\pages  319--330
\endref

\ref 
\key Ia1
\by A\. Iarrobino
\book Compressed algebras: Artin algebras having given socle degrees and maximal length
\bookinfo Trans. Amer. Math. Soc.
\vol 285
\publ A.M.S.
\yr 1984
\pages 337--378
\endref


\ref
\key Ia2
\bysame % A\. Iarrobino
\paper Inverse system of a symbolic power. II. The Waring problem for forms
\jour J. Algebra 
\vol 174 
\yr 1995
\pages  1091--1110
\endref

\ref
\key I--K
\by A\. Iarrobino, V\. Kanev
\book Power sums, Gorenstein algebras, and determinantal loci
\bookinfo L.M.N.
\vol  1721
\publ Springer
\yr 1999
\endref

\ref
\key I--R
\by  A\. Iliev, K\. Ranestad
\paper Canonical curves and varieties of sums of powers of cubic polynomial
\jour J. Algebra
\vol 246
\yr 2001
\pages  385--393
\endref

\ref
\key Mi
\by J.C\. Migliore
\book Introduction to liaison theory and deficiency modules
\bookinfo Progress in Mathematics
\vol 165
\publ Birkh\"auser
\yr 1998
\endref

\ref
\key R--S
\by K\. Ranestad, F.O\. Schreyer
\paper Varieties of sums of powers
\jour J. Reine Angew. Math.
\vol 525
\yr 2000
\pages  147--181
\endref

\ref
\key Ra
\by J\. Rathmann
\paper The uniform position principle for curves in characteristic $p$
\jour Math. Ann.  
\vol 276 
\yr 1987
\pages  565--579
\endref

\ref
\key SD
\by B\. Saint-Donat
\paper On Petri's analysis of the linear system of quadrics through a canonical curve
\jour Math. Ann. 
\vol 206 
\yr 1973
\pages 157--175
\endref

\ref
\key Sch1
\by F.O\. Schreyer
\paper Syzygies of canonical curves and special linear series
\jour Math. Ann.
\vol 275
\yr 1986
\pages  105--137
\endref


\ref
\key Sch2
\bysame % F.O\. Schreyer
\paper A standard basis approach to syzygies of canonical curves
\jour J. Reine Angew. Math.
\vol 421
\yr 1991
\pages  83--123
\endref



\endRefs
\enddocument